\input ./style/arxiv-general.cfg
\documentclass[aap,MSNbibl,seceqn,dvips]{arximspdf}
\makeatletter
   \@ifpackageloaded{graphicx}{}{\usepackage{graphicx}}
\makeatother

\doi{10.1214/14-AAP1076}
\volume{25}
\issue{6}
\pubyear{2015}
\firstpage{3381}
\lastpage{3404}
\docsubty{FLA}

\makeatletter
\newcommand{\rrvert}{\vert}
\newcommand{\llvert}{\vert}
\def\cal{\mathcal}
\newcommand{\eqref}[1]{(\ref{#1})}
\newtheorem{theorem}{Theorem}[section]
\newtheorem{lemma}[theorem]{Lemma}

\newtheorem{propositionn}{Proposition}[section]
\newtheorem{corollary}[theorem]{Corollary}

\def\bD{{\mathbf D}}
\def\bC{{\mathbf C}}
\def\b0{\mathbf0}

\def\sM{{\mathcal M}}

\def\R{\mathbb{R}}
\def\Rp{\mathbb{R}_+}
\def\Z{\mathbb{Z}}
\def\Zp{\mathbb{Z}_+}
\def\N{\mathbb{N}}
\def\E{\mathbb{E}}
\def\P{\mathbb{P}}

\def\sR{{\cal R}}

\def\sZ{{\mathcal Z}}
\def\sZhat{\widehat\sZ}
\def\sZtild{\widetilde\sZ}

\def\sS{{\mathcal S}}

\def\Xhat{\widehat{X}}
\def\Vhat{\widehat{V}}
\def\What{\widehat{W}}
\def\Ehat{\widehat{E}}

\def\Ebar{\overline{E}}
\def\Fbar{\overline{F}}

\def\Qtild{\widetilde{Q}}

\def\wk{\stackrel{w}{\rightarrow}}

\def\Fbar{\overline{F}}
\def\Vbar{\overline{V}}

\makeatother

\begin{document}
\begin{frontmatter}

\title{Diffusion limits for shortest remaining processing time queues
under nonstandard spatial scaling}
\runtitle{Diffusion limits for srpt queues}
\begin{aug}
\author[A]{\fnms{Amber L.}~\snm{Puha}\corref{}\ead[label=e1]{apuha@csusm.edu}\thanksref{T1}}
\runauthor{A. L. Puha}
\thankstext{T1}{Supported in part by ViaSat Inc. and an ROA supplement
to NSF Grant DMS-12-06772.}
\address[A]{Department of Mathematics\\
California State University San Marcos\\
333 S. Twin Oaks Valley Road\\
San Marcos, California 92096-0001\\
USA\\
\printead{e1}}
\affiliation{California State University San Marcos}
\end{aug}

%
\received{\smonth{1} \syear{2014}}
%
\revised{\smonth{7} \syear{2014}}

%
\begin{abstract}
We develop a heavy traffic diffusion limit theorem under nonstandard spatial
scaling for the queue length process in a single server queue employing
shortest remaining processing time (SRPT). For processing time distributions
with unbounded support, it has been shown that standard diffusion
scaling yields
an identically zero limit. We specify an alternative spatial scaling
that produces
a nonzero limit. Our model allows for renewal arrivals and i.i.d.
processing times
satisfying a rapid variation condition. We add a corrective spatial
scale factor to
standard diffusion scaling, and specify conditions under which the
sequence of
unconventionally scaled queue length processes converges in
distribution to
the same nonzero reflected Brownian motion to which the sequence of
conventionally
scaled workload processes converges. Consequently, this corrective spatial
scale factor characterizes the order of magnitude difference between
the queue
length and workload processes of SRPT queues in heavy traffic. It
is determined by the processing time distribution such that the rate at
which it tends to infinity depends on the rate at which the tail of the
processing
time distribution tends to zero. For Weibull processing time
distributions, we restate
this result in a manner that makes the resulting state space collapse
more apparent.
\end{abstract}

%

%
\end{frontmatter}

\section{Introduction}\label{sec1}

We study the heavy traffic behavior of the queue length process in a
shortest remaining processing time (SRPT) queue. We consider a single
server queue with renewal arrivals and independent and identically
distributed processing times that are also independent of the arrival
process. Jobs are served in a nonidling fashion such that at each
instant the job with the shortest remaining processing time is served
at rate one. This is done with preemption so that if a job arrives for
which the total processing time is smaller than that remaining of the
job in service, the job in service is placed on hold and the arriving
job enters service. Therefore, in order to adequately track the state
of the system, it is necessary to keep track of all remaining
processing times of all jobs in the system. We do this using a measure
valued process that at each time has a unit atom at the remaining
processing time of each job in the system. This is introduced formally
in Section~\ref{sec:model}.

Optimality of shortest remaining processing time, in the sense that it
is the queue length minimizer over all nonidling service disciplines,
has been known since the 1960s \cite{sch68,smi76}. One anticipates
that this results at the expense of lengthy delays for jobs with large
total processing times. Hence, sojourn times are naturally of interest
for SRPT. For Markovian
arrivals, the early work Schrage and Miller \cite{schmil66} develops a formula for the
mean response time
in steady state with extended results available in Schassberger \cite{sch90}
and Perera \cite{NEW2}; see Schreiber \cite{sch93} for a survey.
Also see results in Pavlov \cite{NEW1} and Pechinkin \cite{pec86} on steady
state queue length distributions. Recently, Lin, Wierman and Zwart
\cite{ref:LWZ} followed up on the work in \cite{schmil66} by
characterizing the asymptotic behavior of the steady state mean sojourn
time as the traffic intensity approaches one for a large class of
processing time distributions. Interestingly, the rate at which the
mean sojourn time tends to infinity depends on the tail behavior of the
processing time distribution. Results in this spirit were also an
outcome of \cite{ref:sig,ref:DGP}, where a fluid model (formal
functional law of large numbers limit) was proposed and an associate
weak convergence result (functional law of large numbers result) was
stated and proved. There the rate at which a fluid analog of the
sojourn time as a function of the initial processing time tends to
infinity depends on the tail behavior of the processing distribution.
Other somewhat recent studies of SRPT have focused on fairness (e.g.,
\cite{banhar01,wiehar03}) or tail behavior \cite{nun02,nuyzwa06}.

Here we focus on further developing existing diffusion limit results
(functional central limit theorems) for SRPT. The paper \cite{ref:GKP}
contains a diffusion limit theorem for the sequence of measure valued
state descriptors under standard heavy traffic conditions and standard
diffusion scaling. For limiting processing time distributions with
bounded support, the main result in \cite{ref:GKP} indicates that the
limiting measure valued process is a single atom supported at the
supremum of the support of the limiting processing time distribution
for all time. The height of that atom varies randomly in time as
determined by the limiting workload process. More specifically, under
standard heavy traffic conditions, the sequence of conventionally
diffusion scaled workload processes converges in distribution to a
semi-martingale reflected Brownian motion \cite{ref:IW}. The height of
the atom for the measure valued
diffusion limit is then given by the limiting workload process divided
by the supremum of the support of the limiting processing time
distribution. This is analogous to early results for strict priority
queues where in the heavy traffic diffusion limit work piles up in the
lowest priority class \cite{ref:W2}. The main result in~\cite
{ref:GKP} goes on to state that for limiting processing time
distributions with unbounded support, the limiting measure valued
process is identically equal to the zero measure. In particular, the
limiting queue length process is identically equal to zero. Such
behavior had not been observed prior to this for other nonidling
service disciplines. The limiting queue length process is typically
recovered from the limiting workload process via multiplication by a
positive constant, a phenomenon known as state space collapse; see
\cite{ref:G,ref:IW,ref:L,ref:W2}, for instance. The fact that the
sequence of rescaled queue length processes is of lower order magnitude
than the sequence of rescaled workload processes quantifies the extreme
queue length minimizing nature of SRPT.\vadjust{\goodbreak}

A natural follow-up question to the work in \cite{ref:GKP} is whether
or not there is an alternative scaling that can be employed to yield a
nontrivial limit for this unconventionally rescaled queue length
process. If such a limit exists, it would be of interest to describe
how that limit is related to the limiting workload process that arises
under standard diffusion scaling. Here we identify such a nonstandard
scaling for continuous processing time distributions with unbounded
support for which the tails satisfy a rapid variation condition. The
main theorem in this paper, Theorem~\ref{thrm:main}, specifies that
the spatial scaling must be modified by multiplying by a certain
inverse function related to the tails of the first moment of the
processing time distribution. In particular, there is multiplicative
correction factor that must be applied to standard diffusion scaling to
obtain a nontrivial limit. The order of magnitude of that correction
factor depends on the rate at which the tails of the first moment tend
to zero. With this corrective scaling, the limiting process is
identically equal to the limiting workload process that arises under
standard diffusion scaling. Hence, with this corrective scaling factor,
a generalized version of state space collapse holds.

The corrective scaling identified here was inspired by fluid limit
results in \cite{ref:DGP}. The order of magnitude agrees with that of
the left edge of the support of fluid model solutions as time
approaches infinity. This seems to be the first result in the queuing
theory literature where the nature of the scaling depends on the tail
behavior of the processing time distribution. In fact, we only know of
one previous result \cite{wil96} that employs nonstandard scaling. The
scaling in \cite{wil96} is a mixture of conventional fluid (functional
law of large numbers) and conventional diffusion (functional central
limit theorem) scaling.

It is interesting to note that the result in Theorem~\ref{thrm:main}
is consistent with the rapid variation case of
\cite{ref:LWZ}, Theorem~3, as follows. Theorem~3 in \cite{ref:LWZ}
specifies an asymptotic formula for the mean
sojourn time in steady state as the traffic intensity increases to one.
By using the rate at which the traffic intensity
approaches one for standard heavy traffic conditions [see \eqref
{eq:param}], one can informally translate their asymptotic
formula into one indexed by the sequence of systems here. This results
in an asymptotic formula that has the same
order of magnitude as the spatial scaling specified by Theorem~\ref
{thrm:main}. Of course the former is for the steady-state
mean response time, and the latter is for the unconventionally rescaled
queue length process. But, due to Little's law, the queue length and
response time should be of the same order of magnitude.

In general, the inverse function that produces the corrective scaling
is not available in closed form. Hence, the multiplicative constant
contained within it is not immediately available. However, for Weibull
processing time distributions, explicit calculations can be done to
separate the order of magnitude and multiplicative constant. The order
of magnitude is determined by the shape parameter and the
multiplicative constant is given by the scale parameter. This is stated
precisely in Corollary~\ref{cor:Weibull}, which provides an
interesting illustration of the resulting generalize state space collapse.

This raises the next natural question of what happens when the
processing time distributions satisfy a regular variation condition.
The work here does not address that case. The works \cite{ref:DGP} and
\cite{ref:LWZ} suggest that the same function might provide an
appropriate corrective scaling. However, the proof of Theorem~\ref
{thrm:main} does not generalize to that case. The slowly varying nature
of the inverse function that specifies the corrective scaling factor
plays an important role in the proof of Theorem~\ref{thrm:main}; see
\eqref{eq:svrate} and \eqref{eq:rate}. Determining the behavior in
the case of regular variation is work in progress.

In the next section, we precisely define the model and associated
measure valued state descriptor. Then we specify the sequence of
systems and associated asymptotic conditions that they must satisfy.
This allows us to state the main result of the paper, Theorem~\ref
{thrm:main}, and its corollary for Weibull processing time
distributions, Corollary~\ref{cor:Weibull}. The remainder of the paper
contains the proof of the main result.

\subsection{Notation}

Throughout
$\R$ denotes the real numbers, and $\R_+$ denotes the nonnegative
real numbers.
Similarly, $\Z$ denotes the integers, and $\Z_+$ denotes the
nonnegative integers.
Then $\N$ denotes the positive integers. For $a,b\in\R$, $a\wedge b$
and $a\vee b$,
respectively, denote the minimum and maximum of $a$ and $b$. Also, for
$a\in\R$,
$\llvert a\rrvert =(-a)\vee a$ denotes the absolute value of $a$.

We define $\bC(\Rp)$ to be the set of continuous real valued
functions with domain~$\Rp$.
Then $\bC_b(\Rp)$ denotes those elements of $\bC(\Rp)$ that are bounded.
We use the notation $1(\cdot)$ for the function in $\bC_b(\Rp)$ that
is identically equal to one
and $\chi(\cdot)$ for the identity function in $\bC(\Rp)$.

For a Polish space $\sS$, we let $\bD([0,\infty),\sS)$ denote the set
of functions
of time taking values in $\sS$ that are right continuous with finite
left limits.
We endow this space with the Skorohod $J_1$-topology. Then $\bD
([0,\infty)
,\sS)$ is also
a Polish space~\cite{ref:EK}. We denote the function in $\bD
([0,\infty)
,\R)$ that is identically equal
to zero by $0(\cdot)$.

We use the notation $\sM$ for the set of finite, nonnegative Borel
measures on $\Rp$.
The zero measure in $\sM$ is denoted by $\b0$. Traditionally, for
$x\in\Rp$, $\delta_x\in\sM$
is the unit atom at $x$. For $x\in\Rp$, we also define
$\delta_x^+$ to be the measure in $\sM$ that is $\delta_x$ if $x>0$
and $\b0$ otherwise.
Given a Borel measurable function $f\dvtx \Rp\to\R$ and $\zeta\in\sM$,
we let
$\langle f,\zeta\rangle =\int_{\Rp} f(x) \zeta(dx)$, when the
integral exists.
Then $\langle 1,\zeta\rangle $ is the total mass of $\zeta$. We
refer to
$\langle \chi,\zeta\rangle $ as the first moment of $\zeta$. The
set $\sM$
is endowed with the
topology of weak convergence. In particular, for $\{\zeta_n\}_{n\in\N
}\subset\sM$
and $\zeta\in\sM$, $\zeta_n\wk\zeta$ as $n\to\infty$ if and
only if
$\lim_{n\to\infty}\langle g,\zeta_n\rangle =\langle g,\zeta
\rangle $ for all $g\in
\bC_b(\Rp)$.
With this topology, $\sM$ is a Polish space. We denote the function
in $\bD([0,\infty),\sM)$ that is identically equal to the zero measure
by $\b0(\cdot)$.

We use ``$\Rightarrow$'' to denote convergence in distribution of
random elements of a metric space. Following Billingsley \cite{ref:B},
we use $\P$ and $\E$, respectively, to denote the probability
measure and expectation operator associated with whatever space the
relevant random element is defined on. Unless otherwise specified, all
stochastic processes used in this paper are assumed to have paths
that are right continuous with finite left limits (r.c.l.l.).

Finally, following \cite{ref:BGT}, we say that a measurable function
$g\dvtx (0,\infty)\to(0,\infty)$ is rapidly varying of index $\infty$ if
for all
$\varepsilon>0$,
%
\begin{equation}
\label{eq:rvplus} \lim_{x\to\infty}\frac{g((1+\varepsilon)x)}{g(x)}=\infty,
\end{equation}
and is rapidly varying of index $-\infty$ if for all $\varepsilon>0$,
%
\begin{equation}
\label{eq:rvminus} \lim_{x\to\infty}\frac{g((1+\varepsilon)x)}{g(x)}=0.
\end{equation}
Together the functions in these two classes are called \textit{rapidly
varying}.
In addition, $g\dvtx (0,\infty)\to(0,\infty)$ is slowly varying if for
all $c>0$,
%
\begin{equation}
\label{eq:sv} \lim_{x\to\infty}\frac{g(cx)}{g(x)}=1.
\end{equation}

\section{The stochastic model and state descriptor}\label{sec:model}

We consider a $GI/GI/1$ SRPT queue such that the processing time
distribution is continuous, has
unbounded support and the tails satisfy a rapid variation condition. In
particular, jobs arrive
according to a delayed renewal process $E(\cdot)$ with rate $\lambda
\in(0,\infty)$ such that
the interarrival times have finite standard deviation $\sigma_a$ and
$E(0)=0$. Then, for
$t\in[0,\infty)$, $E(t)$ denotes the number of jobs that have arrived to
the system exogenously
by time $t$. Processing times for these jobs
are independent and identically distributed positive random variables
with common continuous
cumulative distribution function $F(\cdot)$, finite mean and finite standard
deviation $\sigma_s\in(0,\infty)$. The sequence $\{v_i\}_{i\in\N}$
of processing times is also
assumed to be independent of the arrival process. For $i\in\N$, the
$i$th job to arrive to the
system has total processing time~$v_i$. For simplicity, we refer to the
$i$th job to arrive to the
system as job $i$, or the $i$th job. We use the notation $v$ to denote
a random variable
that is equal in distribution to a generic processing time. Specifically,
\[
\Fbar(x)=1-F(x)=\P ( v>x ),\qquad x\in\Rp.
\]

Our assumptions include that $\Fbar(\cdot)$ is rapidly varying with
index minus infinity; see \eqref{eq:rvminus}.
We restrict attention to a subset of such processing time distributions
that includes, for example,
Weibull distributions. For this, given $x\in\Rp$, let
\[
S(x)= \frac{1}{{\mathbb E}  [ v 1_{\{ v>x\}}  ]}.
\]
In \cite{ref:DGP}, $s(\cdot)$ is the fluid analog of the
sojourn time of initial jobs as a function of the remaining
processing time at time zero. The notation $S(\cdot)$ is
chosen here to highlight its similarity with $s(\cdot)$
in \cite{ref:DGP} (they differ by factor that tends
to a positive constant as $x$ tends to infinity).
Note that $S(0)=1/{{\mathbb E}  [ v  ]}$. Further,
$S(\cdot)$ is positive, nondecreasing, continuous, unbounded and
rapidly varying with index plus infinity.
Set
%
\begin{equation}
\label{def:Sinv} S^{-1}(y)=\inf\bigl\{ x\in\Rp\dvtx S(x)>y\bigr\}, \qquad y\in\Rp.
\end{equation}
Hence, $S^{-1}(\cdot)$ is positive, strictly increasing,
right continuous, unbounded and slowly varying.
Further, $S(S^{-1}(y))=y$ for all $y\in\Rp$.
We assume that for some $c>1$,
%
\begin{equation}
\label{eq:svrate} \lim_{y\to\infty} \biggl(\frac{S^{-1}( c y )}{S^{-1}(y)}-1 \biggr)\ln
\bigl( S^{-1}(y) \bigr)=0.
\end{equation}
This is an assumption about the rate at which the ratio associated with
the slowly
varying function $S^{-1}(\cdot)$ converges to one. Rate of convergence
conditions such
as this and their implications are discussed more fully in \cite{ref:BGT},
Section~2.3.1.
Here we note that \eqref{eq:svrate} is not satisfied by all slowly
varying functions.
For example, as noted in \cite{ref:BGT}, page 78, \eqref{eq:svrate}
does not hold for
slowly varying functions of the form $\exp (  (\ln(\cdot
) )^{\delta} )$,
where $1/2\le\delta<1$. However, it does hold for many processing
distributions,
including Weibull processing time distributions. That Weibull
processing time
distributions satisfy \eqref{eq:svrate} is demonstrated in Section~\ref{sec:main}.

The reason for assuming \eqref{eq:svrate} is that by \cite{ref:BGT}, Theorem~2.3.3
(originally stated in~\cite{ref:BS}), it follows that for all $\delta
\in\R$,
%
\begin{equation}
\label{eq:rate} \lim_{y\to\infty} \frac{ S^{-1} (  ( S^{-1}(y)
)^{\delta} y ) }{ S^{-1}(y) }=1.
\end{equation}
Recall that $S^{-1}(\cdot)$ is slowly varying so that $\lim_{y\to
\infty}S^{-1}(cy)/S^{-1}(y)=1$ for all $c>0$.
Then \eqref{eq:rate} says that one can replace the constant $c>0$ with
$ ( S^{-1}(y)  )^{\delta}$,
$y\in\Rp$.
As $y$ tends to infinity, this tends to infinity if $\delta>0$ and to
zero if $\delta<0$.
In Section~\ref{sec:main}, \eqref{eq:rate} is used to obtain \eqref
{eq:ratios},
which is in turn used to prove the main theorem of the paper, Theorem~\ref{thrm:main}.

As far as the initial state of the system in concerned, there are
$Q(0)$ jobs in the system at time zero.
Here $Q(0)$ is assumed to be a random variable taking values in $\Zp$.
The time zero remaining processing times for such jobs are the first
$Q(0)$ elements of the sequence
$\{ \tilde v_i\}_{i\in\N}\subset\Rp$. Each member of the sequence
$\{ \tilde v_i\}_{i\in\N}$ is assumed to be a
positive random variable. For $1\le i\le Q(0)$, we refer to the job in
the system at time zero with remaining
processing time $\tilde v_i$ at time zero as initial job $i$, or the
$i$th initial job.
Let $W(0)=\sum_{i=1}^{Q(0)}\tilde v_i$, which is a random variable
taking values in
$\Rp$. Then $W(0)$ corresponds to the total work (measured in units of
processing time) in the system at
time zero. Finally, let $\sZ(0)\in\sM$ be given by
\[
\sZ(0)=\sum_{i=1}^{Q(0)} \delta_{\tilde v_i}^+.
\]
Note that
\[
Q(0)=\bigl\langle 1,\sZ(0)\bigr\rangle\quad \mbox{and}\quad W(0)=\bigl\langle \chi,\sZ (0)
\bigr\rangle.
\]

Jobs are served in a nonidling fashion. In particular, the server does
not idle if there are jobs in the system.
At any given instance at which the system is nonempty, the job with the
shortest remaining processing
time is served at rate one. This is done with preemption so that when a
job arrives to the system that requires less
processing time than that remaining for the job currently in service,
the job in service is placed on hold and the arriving
job enters service immediately. For $1\le i\le Q(0)$ and $t\in
[0,\infty)
$, $\tilde v_i(t)$ denotes the remaining processing time of initial job
$i$ at time $t$. For $1\le i\le E(t)$ and $t\in[0,\infty)$, $v_i(t)$
denotes the remaining processing time of job $i$ at time $t$. So then,
for $t\in[0,\infty)$, let
\[
\sZ(t)=\sum_{i=1}^{Q(0)} \delta_{\tilde v_i(t)}^+
+\sum_{i=1}^{E(t)} \delta_{v_i(t)}^+.
\]
In particular, $\sZ(\cdot)\in\bD([0,\infty),\sM)$ is the associated
measure valued state descriptor. For $t\in[0,\infty)$, let
\[
Q(t)=\bigl\langle 1,\sZ(t)\bigr\rangle\quad \mbox{and}\quad W(t)=\bigl\langle \chi,\sZ (t)
\bigr\rangle.
\]
Then $Q(\cdot)$ and $W(\cdot)$, respectively, denote the queue length
and workload processes.

\section{Statement of the main result}\label{sec:main}

Let $\sR$ be a sequence taking values in $(1,\infty)$ tending to infinity.
Fix a sequence of $GI/GI/1$ SRPT queues indexed by $\sR$ for which the
initial conditions
and stochastic primitive inputs satisfy the conditions specified in
Section~\ref{sec:model}.
We further require that the processing time distributions do not depend
on $r$ and have
common cumulative distribution function $F(\cdot)$. We place a
superscript $r$ on all parameters and processes associated with the
$r$th system. So
then for each $r\in\sR$, we have $\lambda^r$, $\sigma_a^r$,
$E^r(\cdot)$, $\sZ^r(\cdot)$,
$Q^r(\cdot)$ and $W^r(\cdot)$, which may depend on $r$, but
$F^r(\cdot)=F(\cdot)$ for
all $r\in\sR$. Also, for $r\in\sR$, set
\[
\rho^r=\lambda^r{\mathbb E} [v ].
\]
For convenience later on, for $r\in\sR$ and $x\in\Rp$, we also define
\[
\rho_x^r=\lambda^r \E [ v 1_{\{ v\le x\}}
].
\]
Then, for $r\in\sR$ and $x\in\Rp$,
%
\begin{equation}
\label{eq:rhominusrhox} \rho^r-\rho_x^r=
\frac{\lambda^r}{S(x)}.
\end{equation}

We assume that the stochastic primitive inputs satisfy the following
asymptotic heavy traffic conditions.
For some $\kappa\in\R$, as $r\to\infty$,
%
\begin{equation}
\label{eq:param} \sigma_a^r\to\sigma_a\quad
\mbox{and}\quad r \bigl(\rho^r-1 \bigr)\to\kappa.
\end{equation}
Then it follows that $\lambda^r\to\lambda$ as $r\to\infty$, where
$\lambda=1/\E[v]$.
For $r\in\sR$ and $t\in[0,\infty)$, let
\[
\Ebar^r(t)=\frac{E^r(r^2 t)}{r^2} \quad\mbox{and}\quad \Ehat^r(t)=
\frac{ E^r(r^2 t)-\lambda^rr^2t}{r}.
\]
Also assume that as $r\to\infty$,
%
\begin{equation}
\label{eq:FCLTE} \Ehat^r(\cdot)\Rightarrow E^*(\cdot),
\end{equation}
where $E^*(\cdot)$ is a Brownian motion starting from zero
with drift zero and variance $(\lambda)^3(\sigma_a)^2$.
This implies a functional weak law of large numbers for the
exogenous arrival process. Specifically, set
$\lambda(t)=\lambda t$ for $t\in[0,\infty)$. Then, as $r\to\infty$,
%
\begin{equation}
\label{eq:FLLNE} \Ebar^r(\cdot)\Rightarrow\lambda(\cdot).
\end{equation}

Given $r\in\sR$, let
%
\begin{equation}
\label{def:cr} c^r=S^{-1} (r ).
\end{equation}
For $r\in\sR$ and $t\in[0,\infty)$, set
\[
\Qtild^r(t)=\frac{ c^r Q^r(r^2t)}{r}\quad \mbox{and} \quad\sZtild^r(t)=
\frac{ c^r \sZ^r(r^2 t)}{r}.
\]
Also, for $r\in\sR$ and $t\in[0,\infty)$, set
\[
\What^r(t)=\frac{W^r(r^2t)}{r} \quad\mbox{and} \quad \sZhat^r(t)=
\frac{\sZ^r(r^2 t)}{r}.
\]
Then, the ``hat'' notation corresponds to processes under standard
diffusion scaling
and the ``tilde'' notation corresponds to processes under the
nonstandard scaling
consisting of standard diffusion scaling multiplied by the spatial
correction factor
$c^r$, $r\in\sR$. Note that $\lim_{r\to\infty}c^r=\infty$.
Assume that for some random variable $W_0$ that is finite almost surely,
as $r\to\infty$,
%
\begin{equation}
\label{eq:initial1} \bigl(\What^r(0),\Qtild^r(0) \bigr)
\Rightarrow (W_0, W_0 ).
\end{equation}
For $r\in\sR$ and $\varepsilon>0$, let
%
\begin{equation}
\label{def:lu} l_{\varepsilon}^r=S^{-1} \bigl(r
\bigl(c^r\bigr)^{-2-\varepsilon} \bigr) \quad\mbox{and}\quad u_{\varepsilon}^r=S^{-1}
\bigl(r \bigl(c^r\bigr)^{2+\varepsilon} \bigr).
\end{equation}
Then, for $\varepsilon>0$ and $r\in\sR$, we have that
$0<l_{\varepsilon}^r<c^r<u_{\varepsilon}^r<\infty$.
Also, for all $\varepsilon>0$, $\lim_{r\to\infty}l_{\varepsilon
}^r=\lim_{r\to\infty}u_{\varepsilon}^r=\infty$.
Further, by \eqref{eq:rate}, \eqref{def:cr} and \eqref{def:lu}, for
each $\varepsilon>0$,
%
\begin{equation}
\label{eq:ratios} \lim_{r\to\infty} \frac{c^r}{l_{\varepsilon}^r}=1 \quad\mbox{and}\quad \lim
_{r\to\infty} \frac{c^r}{u_{\varepsilon}^r}=1.
\end{equation}
The proof of the main result (Theorem~\ref{thrm:main}) will proceed by
demonstrating for any
given $\varepsilon>0$, the contribution to the total mass under the
unconventional scaling
and to the work under the conventional scaling asymptotically
concentrates in $(l_{\varepsilon}^r,u_{\varepsilon}^r]$ as $r\to
\infty$. Therefore, we further assume that
for all $\varepsilon>0$, as $r\to\infty$,
%
\begin{equation}
\label{eq:initial2} \bigl\langle (1\vee\chi ) 1_{[0,l_{\varepsilon}^r]},\sZtild
^r(0)\bigr\rangle \Rightarrow0 \quad\mbox{and}\quad \bigl\langle
\chi1_{(u_{\varepsilon}^r,\infty)},\sZhat^r(0)\bigr\rangle \Rightarrow0.
\end{equation}

\begin{theorem}\label{thrm:main}
Assume that \eqref{eq:param}, \eqref{eq:FCLTE}, \eqref{eq:initial1}
and \eqref{eq:initial2} hold. As $r\to\infty$,
\[
\bigl(\Qtild^r(\cdot),\What^r(\cdot) \bigr)\Rightarrow
\bigl(W^*(\cdot ),W^*(\cdot)\bigr),
\]
where $W^*(\cdot)$ is a reflected Brownian motion with drift $\kappa$
and variance\break 
$\lambda((\sigma_a)^2+(\sigma_s)^2)$ such that $W^*(0)$ is equal in
distribution to $W_0$.
\end{theorem}

For the class of processing time distributions that satisfy the rapid
variation condition \eqref{eq:svrate},
Theorem~\ref{thrm:main} implies that the asymptotic order of magnitude
difference between the $\sR$
indexed queue length and workload processes in heavy traffic is given
by $c^r=S^{-1}(r)$,
$r\in\sR$. Through \eqref{def:Sinv}, the order of magnitude of the
correction factor $c^r$, $r\in\sR$, is
determined by the rate at which the tail of the first moment of the
processing time distribution tends to zero.

One can view Theorem~\ref{thrm:main} as a generalized state space
collapse result with a
multiplicative lifting factor of one; that is, the heavy traffic limit
of the unconventionally rescaled queue length process
is one times the heavy traffic limit of the conventionally rescaled
workload process. The proof of Theorem~\ref{thrm:main}
given in Section~\ref{sec:proof} provides insight into how this
phenomenon manifests itself. We give an informal overview
there as well. Another way to view this result is that the sequence of
spatial correction factors $\{c^r\}_{r\in\sR}$, has
embedded in it both the order of magnitude difference between the $\sR
$ indexed queue length and
workload processes in heavy traffic and the reciprocal of the
multiplicative lifting map. For many processing
time distributions that are of interest in practice, one can compute
these explicitly. We illustrate this in the following
corollary.

In the following corollary, we consider Weibull processing time
distributions with positive
shape parameter $\alpha$ and positive rate parameter $\beta$. For
these processing time distributions, the corollary
precisely identifies the order of magnitude of the corrective spatial
scaling factor as $\sqrt[\alpha]{\ln r}$. It also identifies
what can be viewed as a state space collapse lifting map that obtains
the limit of the sequence of diffusion scaled queue
length processes with the $r$th member multiplied by $\sqrt[\alpha
]{\ln r}$ from the limit of the sequence of diffusion scaled workload
process via multiplication by the rate parameter $\beta$. In this
regard, it is interesting to note that multiplication of the limiting
workload process by $\beta$ is not the
same as division by the mean processing time, except in the exponential
case $\alpha=1$. Indeed, the mean processing
time is given by $\Gamma(1+\alpha)/\beta$, where $\Gamma(t)=\int_{\Rp} x^{t-1}\exp(-x)\,dx$, $t\in(0,\infty)$, denotes
the gamma function. Note that $\Gamma(1+\alpha)<1$ for $0<\alpha<1$
and $\Gamma(1+\alpha)>1$ for $\alpha>1$.
Then, under this nonstandard spatial scaling, the limiting residual
processing time per job in the system $1/\beta$
exceeds the mean processing time for $0<\alpha<1$. The opposite is
true for $\alpha>1$.

\begin{corollary}\label{cor:Weibull} Let $\alpha,\beta>0$.
Assume that \eqref{eq:param}, \eqref{eq:FCLTE}, \eqref{eq:initial1}
and \eqref{eq:initial2} hold and that
$\Fbar(x)=\exp(-(\beta x)^{-\alpha})$, $x\in\Rp$ (so that the
processing time distribution is Weibull distributed with rate parameter
$\beta>0$ and shape parameter $\alpha>0$). Then, as $r\to\infty$,
\[
\frac{ \sqrt[\alpha]{\ln(r)} Q^r(r^2\cdot)}{r} \Rightarrow\beta W^*(\cdot),
\]
where $W^*(\cdot)$ is a reflected Brownian motion with drift $\kappa$
and variance\break $\lambda((\sigma_a)^2+(\sigma_s)^2)$ such that $W^*(0)$
is equal in distribution to $W_0$.
\end{corollary}

\begin{pf}
Fix $\alpha,\beta>0$. We begin by more precisely determining the
asymptotic behavior of $S^{-1}(\cdot)$; see \eqref{eq:Weibull} below.
Then we use this asymptotic behavior to verify \eqref{eq:svrate} so
that we may apply
Theorem~\ref{thrm:main}.\vspace*{1pt} The continuous mapping theorem together with
\eqref{eq:Weibull}, then allows us to
replace $c^r=S^{-1}(r)$ with $\sqrt[\alpha]{\ln r}/\beta$ and then
to multiply by the constant $\beta$ to
obtain the desired conclusion.

For $x\in\Rp$,
\[
\frac{1}{S(x)}={\mathbb E} [ v 1_{\{ v>x\}} ]=x\Fbar (x)+\int
_x^{\infty} \Fbar(y) \,dy \ge \frac{x}{\exp ( (\beta x)^{\alpha}  )}.
\]
Using L'Hopital's rule, one can verify that
\[
\lim_{x\to\infty}\frac{\exp((\beta x)^{\alpha})}{x S(x)}= \lim_{x\to\infty}
\frac{{\mathbb E}  [ v 1_{\{ v>x\}}
]}{x\Fbar(x)}=1.
\]
Fix $\delta\in(0,1)$. Then there exists $X\in\Rp$ such that for all $x>X$,
\begin{eqnarray*}
(1-\delta)\exp \bigl( \bigl( (1-\delta)\beta x \bigr)^{\alpha} \bigr) &
\le& \frac{(1-\delta) \exp ( (\beta x)^{\alpha}  )}{x}
\\
&\le& S(x)
\\
&\le& \frac{\exp ( (\beta x)^{\alpha} )}{x}
\\
&\le& \exp \bigl( (\beta x)^{\alpha} \bigr).
\end{eqnarray*}
So then it follows that there exists $Y\in\Rp$ such that for $y>Y$,
\[
\frac{\sqrt[\alpha]{\ln(y)}}{\beta}\le S^{-1}(y)\le\frac{\sqrt
[\alpha]{\ln ({y}/{(1-\delta)} )}}{(1-\delta)\beta}.
\]
Hence
%
\begin{equation}
\label{eq:Weibull} \lim_{y\to\infty} \frac{ \beta S^{-1}(y)}{ \sqrt[\alpha]{\ln y}}=1.
\end{equation}

Fix $c>1$. For $y>1$, we have
\begin{eqnarray*}
&& \biggl( \frac{\sqrt[\alpha]{\ln(cy)}}{\sqrt[\alpha]{\ln(y)}} -1 \biggr)\ln \biggl(\frac{\sqrt[\alpha]{\ln(y)}}{\beta} \biggr)
\\
&&\qquad = \biggl( \sqrt[\alpha]{1+\frac{\ln(c)}{\ln(y)}} -\sqrt[\alpha ]{1} \biggr)
\biggl( \frac{\ln (\ln(y) )}{\alpha}-\ln \beta \biggr).
\end{eqnarray*}
Set $h(z)=\sqrt[\alpha]{1+z}$, $z\in(-1,\infty)$. Using Taylor's
remainder theorem
and the fact that $h'(\cdot)$ is continuous in a neighborhood of the
origin, there exists $B,\delta>0$ such
that for all $\llvert  z\rrvert  <\delta$,
\[
1-B\llvert z\rrvert \le h(z)\le1+B\llvert z\rrvert .
\]
So then for all $y$ sufficiently larger than 1,
\[
0\le \biggl(\sqrt[\alpha]{1+\frac{\ln(c)}{\ln(y)}} -\sqrt[\alpha ]{1} \biggr)\ln
\biggl(\frac{\sqrt[\alpha]{\ln(y)}}{\beta} \biggr) \le \frac{B \ln(c) }{\ln(y)} \biggl(
\frac{\ln (\ln(y)
)}{\alpha}-\ln\beta \biggr).
\]
Hence
%
\begin{equation}
\label{eq:svrateWeibull} \lim_{y\to\infty} \biggl( \frac{\sqrt[\alpha]{\ln(cy)}}{\sqrt
[\alpha]{\ln(y)}} -1 \biggr)
\ln \biggl(\frac{\sqrt[\alpha]{\ln
(y)}}{\beta} \biggr)=0.
\end{equation}
Combining \eqref{eq:Weibull} and \eqref{eq:svrateWeibull} implies
\eqref{eq:svrate}
for $S^{-1}(\cdot)$. Hence the result follows from Thoerem \ref
{thrm:main}, \eqref{eq:Weibull} and the continuous
mapping theorem.
\end{pf}

\section{Proof of Theorem \texorpdfstring{\protect\ref{thrm:main}}{3.1}}\label{sec:proof}

Here we state the main facts that will be proved in subsequent sections
in order to verify Theorem~\ref{thrm:main}. Then we prove Theorem~\ref
{thrm:main}
using these facts.

Henceforth, we assume that we have a sequence of $\sR$ indexed $GI/GI/1$
SRPT queues satisfying
the conditions in Section~\ref{sec:main} and that $W^*(\cdot)$
denotes a semi-martingale reflected
Brownian motion with drift $\kappa$ and variance $\lambda((\sigma
_a)^2+(\sigma_s)^2)$ such that
$W^*(0)$ is equal in distribution to $W_0$. Then, by \cite{ref:IW}, as
$r\to\infty$,
%
\begin{equation}
\label{eq:W} \widehat W^r(\cdot)\Rightarrow W^*(\cdot).
\end{equation}
In Section~\ref{sec:SRPTW}, we state and prove Lemma~\ref{lem:WorkAbove}.
This together with the fact that $c^r< u_{\varepsilon}^r$ for all
$\varepsilon>0$
and $r\in\sR$ implies that for all $\varepsilon>0$, as $r\to\infty$,
%
\begin{equation}
\label{eq:AboveU} \bigl\langle 1_{(u_{\varepsilon}^r,\infty)},\sZtild^r(\cdot)\bigr
\rangle \Rightarrow 0(\cdot) \quad\mbox{and}\quad \bigl\langle \chi1_{(u_{\varepsilon}^r,\infty)},
\sZhat^r(\cdot)\bigr\rangle \Rightarrow0(\cdot).
\end{equation}
In Section~\ref{sec:SRPTQ}, we state and prove Lemma~\ref{lem:BelowL}.
This implies that for all $\varepsilon>0$, as $r\to\infty$,
%
\begin{equation}
\label{eq:BelowL} \bigl\langle 1_{[0,l_{\varepsilon}^r]},\sZtild^r(\cdot)\bigr
\rangle \Rightarrow 0(\cdot) \quad\mbox{and}\quad \bigl\langle \chi1_{[0,l_{\varepsilon}^r]},
\sZhat^r(\cdot)\bigr\rangle \Rightarrow 0(\cdot).
\end{equation}
The asymptotic behavior summarized in \eqref{eq:AboveU} and \eqref{eq:BelowL}
is used below in the proof of Theorem~\ref{thrm:main}.

Before proceeding to prove Theorem~\ref{thrm:main}, we provide an
overview, which
provides some insight into how the state space collapse that it implies
arises. For this, let $\varepsilon>0$.
Then \eqref{eq:AboveU} and \eqref{eq:BelowL} imply that in heavy
traffic the jobs that contribute to the
unconventionally rescaled queue length process or to the conventionally
rescaled workload process have
residual processing times that asymptotically concentrate in
$(l_{\varepsilon}^r,u_{\varepsilon}^r]$ as $r\to\infty$.
For each $r\in\sR$, this interval contains the scale factor $c^r$.
The interval itself is shifting out to infinity as
$r\to\infty$. However, since the workload process converges to a
nondegenerate limit under diffusion scaling,
the number of jobs with residual service time in this interval must
tend to zero on diffusion scale. That the diffusion
scaled queue length has a zero limit was shown rigorously in \cite
{ref:GKP}, which implies that the diffusion
scaled measure valued state descriptor has a zero limit as well.
However, due to \eqref{eq:ratios}, all members
of this interval are of order $c^r$. In particular, each job with
residual processing time in this interval contributes
order $c^r$ to the diffusion scaled workload process. So then, since
jobs with residual service time outside of
$(l_{\varepsilon}^r,u_{\varepsilon}^r]$ do not asymptotically
contribute to the unconventionally rescaled queue length process,
it should follow that as $r\to\infty$
\[
c^r \widehat{Q}^r(\cdot)\approx\widehat{W}^r(
\cdot).
\]
The proof of Theorem~\ref{thrm:main} given next demonstrates this in
precise terms, and thereby validates this
line of reasoning.

\begin{pf*}{Proof of Theorem~\ref{thrm:main}}
We have that for all $\varepsilon>0$, $r\in\sR$ and $t\in[0,\infty)$,
\[
l_{\varepsilon}^r \bigl\langle 1_{(l_{\varepsilon}^r,u_{\varepsilon
}^r]},
\sZhat^r(t)\bigr\rangle \le \bigl\langle \chi1_{(l_{\varepsilon}^r,u_{\varepsilon}^r]},\sZhat
^r(t)\bigr\rangle \le u_{\varepsilon}^r\bigl\langle
1_{(l_{\varepsilon}^r,u_{\varepsilon
}^r]},\sZhat^r(t)\bigr\rangle.
\]
Then, for all $\varepsilon>0$, $r\in\sR$ and $t\in[0,\infty)$,
%
\begin{equation}
\label{eq:Squeeze} \frac{c^r }{ u_{\varepsilon}^r} \bigl\langle \chi1_{(l_{\varepsilon
}^r,u_{\varepsilon}^r]},
\sZhat^r(t)\bigr\rangle \le \bigl\langle 1_{(l_{\varepsilon}^r,u_{\varepsilon}^r]},\sZtild
^r(t)\bigr\rangle \le \frac{c^r }{l_{\varepsilon}^r} \bigl\langle
\chi1_{(l_{\varepsilon
}^r,u_{\varepsilon}^r]},\sZhat^r(t)\bigr\rangle.
\end{equation}
Fix $T,\varepsilon, \eta, \delta>0$. Given $r\in\sR$, let
\begin{eqnarray*}
\Omega_1^r&=& \Bigl\{ \sup_{t\in[0,T]}
\bigl\langle 1_{(u_{\varepsilon
}^r,\infty)},\sZtild^r(t)\bigr\rangle <\delta/3 \Bigr
\}\cap \Bigl\{\sup_{t\in[0,T]}\bigl\langle \chi1_{(u_{\varepsilon}^r,\infty
)},
\sZhat^r(t)\bigr\rangle <\delta/3 \Bigr\},
\\
\Omega_2^r&=& \Bigl\{ \sup_{t\in[0,T]} \bigl
\langle 1_{[0,l_{\varepsilon
}^r]},\sZtild^r(t)\bigr\rangle <\delta/3 \Bigr\}
\cap \Bigl\{\sup_{t\in[0,T]}\bigl\langle \chi1_{[0,l_{\varepsilon
}^r]},\sZhat
^r(t)\bigr\rangle <\delta/3 \Bigr\}.
\end{eqnarray*}
By \eqref{eq:AboveU} and \eqref{eq:BelowL},
%
\begin{equation}
\label{eq:ProbOne} \lim_{r\to\infty}\P \bigl( \Omega_1^r
\cap\Omega_2^r \bigr)=1.
\end{equation}
By \eqref{eq:Squeeze}, for each $r\in\sR$, on $\Omega_1^r\cap
\Omega_2^r$, for all $t\in[0,T]$,
\[
\frac{c^r }{ u_{\varepsilon}^r} \What^r(t)-\frac{2\delta}{3} \le
\Qtild^r(t) \le \frac{c^r }{l_{\varepsilon}^r}\What^r(t) +
\frac{2\delta}{3}.
\]
Then, for each $r\in\sR$, on $\Omega_1^r\cap\Omega_2^r$, for all
$t\in[0,T]$,
%
\begin{equation}
\label{eq:Squeeze2} \biggl(\frac{c^r }{ u_{\varepsilon}^r}-1 \biggr)\What^r(t)-
\frac
{2\delta}{3} \le \Qtild^r(t)-\What^r(t) \le \biggl(
\frac{c^r }{l_{\varepsilon}^r}-1 \biggr) \What^r(t)+\frac
{2\delta}{3}.
\end{equation}

Given $r\in\sR$ and $M\in\N$, let
\[
\Omega^r(M)= \Bigl\{ \sup_{t\in[0,T]}
\What^r(t)<M \Bigr\} \quad\mbox{and}\quad \Omega(M)= \Bigl\{\sup
_{t\in[0,T]} W^*(t)<M \Bigr\}.
\]
Since $W^*(\cdot)$ is continuous almost surely,
\[
\P \biggl( \bigcup_{M\in\N} \Omega(M) \biggr)=1.
\]
Hence, there exists $M_{\eta}\in\N$ such that
\[
\P \bigl( \Omega(M_{\eta}) \bigr)\ge1-\eta.
\]
Then, by \eqref{eq:W} and the Portmanteau theorem,
\[
\liminf_{r\to\infty} \P \bigl( \Omega^r(M_{\eta})
\bigr)\ge 1-\eta.
\]
This together with \eqref{eq:ProbOne} implies that
%
\begin{equation}
\label{eq:ProbEta} \liminf_{r\to\infty} \P \bigl( \Omega_1^r
\cap\Omega_2^r\cap \Omega^r(M_{\eta})
\bigr)\ge1-\eta.
\end{equation}
Further, by \eqref{eq:Squeeze2}, for each $r\in\sR$, on $\Omega
_1^r\cap\Omega_2^r\cap\Omega^r(M_{\eta})$,
\[
\sup_{t\in[0,T]}\bigl\llvert \Qtild^r(t)-
\What^r(t)\bigr\rrvert \le \max \biggl( \frac{c^r }{l_{\varepsilon}^r}-1, 1-
\frac{c^r }{
u_{\varepsilon}^r} \biggr)M_{\eta}+\frac{2\delta}{3}.
\]
By \eqref{eq:ratios}, there exists $R\in\sR$ such that for all $r>R$,
\[
\max \biggl( \frac{c^r }{l_{\varepsilon}^r}-1, 1-\frac{c^r }{
u_{\varepsilon}^r} \biggr)\le
\frac{\delta}{3M_{\eta}}.
\]
Then, for each $r>R$, on $\Omega_1^r\cap\Omega_2^r\cap\Omega
^r(M_{\eta})$,
\[
\sup_{t\in[0,T]}\bigl\llvert \Qtild^r(t)-
\What^r(t)\bigr\rrvert \le\delta.
\]
Hence, by \eqref{eq:ProbEta},
\[
\liminf_{r\to\infty}\P \Bigl(\sup_{t\in[0,T]}\bigl\llvert
\Qtild ^r(t)-\What^r(t)\bigr\rrvert \le\delta \Bigr)\ge1-
\eta.
\]
Since $T,\eta,\delta>0$ were arbitrary,
\[
\Qtild^r(\cdot)-\What^r(\cdot)\Rightarrow0(\cdot).
\]
This together with \eqref{eq:W} and the converging together lemma
completes the proof.
\end{pf*}

\section{Verification of \texorpdfstring{\protect\eqref{eq:AboveU}}{(4.2)} and \texorpdfstring{\protect\eqref{eq:BelowL}}{(4.3)}}

Theorem~\ref{thrm:main} was proved in Section~\ref{sec:proof} as a
consequence of
\eqref{eq:AboveU} and \eqref{eq:BelowL} and other facts already
established in the paper.
The remainder of the paper is devoted to stating and proving the two
lemmas that imply
\eqref{eq:AboveU} and \eqref{eq:BelowL}, namely Lemmas \ref
{lem:WorkAbove} and \ref{lem:BelowL}.

\subsection{Workload process tail behavior}\label{sec:SRPTW}

In this section we prove Lemma~\ref{lem:WorkAbove}, which implies
\eqref{eq:AboveU}.
The tail behavior asserted here is relatively easy to verify since it
is simply
a manifestation of the scaling. This is evident in the proof given below.

%

\begin{lemma}\label{lem:WorkAbove}
For all $\varepsilon>0$, as $r\to\infty$,
%
\begin{equation}
\label{eq:ubzero} \bigl\langle \chi1_{(u_{\varepsilon}^r,\infty)},\sZhat^r(\cdot)\bigr
\rangle \Rightarrow0(\cdot).
\end{equation}
\end{lemma}

\begin{pf} Fix $\varepsilon>0$. For $r\in\sR$ and $t\in[0,\infty)$,
$\tilde v_i^r(t)\le\tilde v_i^r$ for all $1\le i\le Q^r(0)$
and $v_i^r(t)\le v_i$ for all $1\le i\le E^r(r^2t)$. Hence, for $r\in
\sR$,
%
\begin{equation}
\label{eq:ub} \bigl\langle \chi1_{(u_{\varepsilon}^r,\infty)},\sZhat^r(\cdot)\bigr
\rangle \le \bigl\langle \chi1_{(u_{\varepsilon}^r,\infty)},\sZhat^r(0)\bigr\rangle
+ \frac{1}{r}\sum_{i=1}^{r^2\Ebar^r(\cdot)}
v_i 1_{\{ v_i >
u_{\varepsilon}^r\} }.
\end{equation}
%
Further, for $r\in\sR$,
\[
\frac{1}{r}\sum_{i=1}^{r^2\Ebar^r(\cdot)}
v_i 1_{\{ v_i>
u_{\varepsilon}^r\} } = \frac{1}{r} \Biggl(\sum
_{i=1}^{r^2\Ebar^r(\cdot)} v_i 1_{\{ v_i>
u_{\varepsilon}^r\} } -
r^2\lambda^r(\cdot){\mathbb E} [ v 1_{\{ v > u_{\varepsilon
}^r\} } ]
\Biggr) +\frac{r\lambda^r(\cdot)}{S(u_{\varepsilon}^r)}.
\]
By \eqref{def:lu}, \eqref{eq:param} and $\lim_{r\to\infty}
c^r=\infty$,
%
\begin{equation}
\label{eq:zero} \lim_{r\to\infty}\frac{r\lambda^r(\cdot)}{S(u_{\varepsilon
}^r)}=\lim
_{r\to\infty}\frac{\lambda^r(\cdot
)}{(c^r)^{2+\varepsilon}}=0(\cdot).
\end{equation}
Further, as $r$ tends to infinity, $\E [v1_{\{v>u_{\varepsilon
}^r\} } ]$ and ${\mathbb E} [ v^21_{\{v>u_{\varepsilon}^r\}
}  ]$ converge to zero since $\lim_{r\to\infty}u_{\varepsilon
}^r=\infty$. Hence, by Proposition~\ref{prop:basic}, as $r\to\infty$,
\[
\frac{1}{r} \Biggl(\sum_{i=1}^{r^2\Ebar^r(\cdot)}
v_i 1_{\{ v_i>
u_{\varepsilon}^r\} } - r^2\lambda^r(\cdot){
\mathbb E} [ v1_{\{
v>u_{\varepsilon}^r\} } ] \Biggr) \Rightarrow0(\cdot).
\]
Therefore, as $r\to\infty$,
\[
\frac{1}{r}\sum_{i=1}^{r^2\Ebar^r(\cdot)}
v_i 1_{\{ v_i>
u_{\varepsilon}^r\} } \Rightarrow0(\cdot).
\]
Combining this with \eqref{eq:initial2} and \eqref{eq:ub} implies
\eqref{eq:ubzero}.
\end{pf}

\subsection{Behavior in large neighborhoods of the origin}

In this section, we prove the following lemma, which implies \eqref{eq:BelowL}.

\begin{lemma}\label{lem:BelowL}
For all $\varepsilon>0$, as $r\to\infty$,
%
\begin{equation}
\label{eq:lbzero} \bigl\langle ( 1\vee\chi ) 1_{[0,l_{\varepsilon}^r]},\sZtild
^r(\cdot)\bigr\rangle \Rightarrow0(\cdot).
\end{equation}
\end{lemma}

The behavior asserted in Lemma~\ref{lem:BelowL} is more subtle than
that asserted
in Lemma~\ref{lem:WorkAbove} since it relies on the SRPT processing dynamics.
Key elements used in verifying this result are asymptotics obtained for
the duration of busy periods for large neighborhoods of the origin; see
Lemmas \ref{lem:ThetaX}
and \ref{lem:theta}. Such results are refinements of \cite{ref:GKP},
(4.9), where the neighborhood
of the origin does not grow with $r\in\sR$, and a slower rate of
convergence to zero is verified for
fixed width neighborhoods of the origin. Equations \eqref
{eq:RepWorkBelowX1} and
\eqref{eq:RepWorkBelowX2} developed below play a central role in
proving these rate of convergence
results. They exploit the nonidling nature of SRPT as well as the order
in which jobs are processed.

Once Lemmas \ref{lem:ThetaX} and \ref{lem:theta} are established, we verify
that the total mass in a fixed width neighborhood of the origin
converges to zero; see Lemma~\ref{lem:Queuex}. The proof of Lemma~\ref
{lem:Queuex} utilizes
an inequality similar in spirit to \eqref{eq:RepWorkBelowX1}, but for
total mass
rather than the total amount of work; see \eqref{eq:QueueRep}. This
inequality is less precise than \eqref{eq:RepWorkBelowX1}
since knowing how many time units the server has spent processing work
does not exactly prescribe
the number of jobs that exit the system during that timeframe. However,
by fixing the
width of the neighborhood of the origin, one can utilize this dynamic
inequality together with
the result in Lemma~\ref{lem:ThetaX} to obtain the desired conclusion.

The final step is to verify that the total amount of work in a growing
neighborhood of the origin
tends to zero; see Lemma~\ref{lem:Worktild}. For this, we return to
\eqref{eq:RepWorkBelowX2}
multiplied by the corrective spatial scaling factor $c^r$ and with $x$
taken to be $l_{\varepsilon}^r$, $\varepsilon>0$
and $r\in\sR$. This yields an upper bound on the desired quantity.
Then we need to verify that all
terms on the right-hand side tend to zero. In particular, we must
verify that the net change over
certain busy periods of what could be referred to as centered truncated
load processes tends to zero
sufficiently fast. This is addressed by Lemma~\ref{lem:Holder}. Since
these centered, truncated
load processes converge to Brownian motion (as noted in the \hyperref[app]{Appendix}),
the proof strategy is to use H\"older
continuity of Brownian motion to bound such differences by quantities
involving the duration of the busy period.
This allows one to utilize the asymptotics obtained in Lemma~\ref
{lem:theta} to prove Lemma~\ref{lem:Holder}.
The result in Lemma~\ref{lem:Holder} is combined with other facts in
order to prove Lemma~\ref{lem:Worktild}
at the end of Section~\ref{sec:SRPTWBelow}.

For completeness, we write out the proof of Lemma~\ref{lem:BelowL} as
a consequence of Lemmas \ref{lem:Queuex}
and \ref{lem:Worktild} here.

\begin{pf*}{Proof of Lemma~\ref{lem:BelowL}}
Fix $\varepsilon>0$.
Then, given $r\in\sR$,
\[
\bigl\langle ( 1\vee\chi ) 1_{[0,l_{\varepsilon}^r]}, \sZtild ^r(\cdot) \bigr
\rangle \le \bigl\langle 1_{[0,1]}, \sZtild^r(\cdot) \bigr
\rangle + \bigl\langle \chi1_{[0,l_{\varepsilon}^r]}, \sZtild^r(\cdot) \bigr
\rangle.
\]
This together with Lemmas \ref{lem:Queuex} and \ref{lem:Worktild}
immediately implies \eqref{eq:lbzero}.
\end{pf*}

The remainder of this section contains the statements and proofs of
Lemmas \ref{lem:Queuex} and \ref{lem:Worktild}.

\subsubsection{Asymptotics for busy period durations}\label{sec:Busy}

For $x\in\Rp$, $r\in\sR$ and $t\in[0,\infty)$, let
\[
\tau^r(t,x)=\sup\bigl\{ s\in[0,t] \dvtx \bigl\langle
1_{[0,x]},\sZtild ^r(s)\bigr\rangle =0\bigr\} \quad\mbox{and}\quad
\theta^r(t,x)=t-\tau^r(t,x).
\]
Given $x\in\Rp$, $r\in\sR$ and $t\in[0,\infty)$, $\theta^r(t,x)$
represents the amount of time
that has elapsed since the $r$th system had no jobs with residual
processing time in $[0,x]$.
In particular, given $x\in\Rp$, $r\in\sR$ and $t\in[0,\infty)$,
$\langle
\chi1_{[0,x]},\break \sZ^r(r^2s)\rangle >0$
for all $s\in(\tau^r(t,x),t]$. Hence, during the time interval
$(r^2\tau^r(t,x),\break r^2t]$
the server in the $r$th system is busy and devoted to serving jobs with
remaining processing
time in $[0,x]$. Hence, for each $x\in\Rp$, $r\in\sR$ and $t\in
[0,\infty)$,
\begin{eqnarray*}
\bigl\langle \chi1_{[0,x]},\sZ^r\bigl(r^2t\bigr)
\bigr\rangle &=& \bigl\langle \chi1_{[0,x]},\sZ^r
\bigl(r^2\tau^r(t,x)\bigr)\bigr\rangle \\
&&{}+ \sum
_{i=
E^r(r^2\tau
^r(t,x))+1}^{E^r(r^2t)}v_i 1_{\{v_i\le x \}}
-r^2 \theta^r(t,x).
\end{eqnarray*}
For $x\in\Rp$, $r\in\sR$ and $t\in[0,\infty)$, set
\begin{eqnarray*}
V_x^r(t)&=&\sum_{i=1}^{E^r(t)}v_i1_{\{v_i\le x\}},
\\
\Vbar_x^r(t)&=&\frac{V_x^r(r^2t)}{r^2},
\\
\Vhat_x^r(t)&=&\frac{1}{r} \bigl(
V_x^r\bigl(r^2t\bigr) -
\rho_x^rr^2 t \bigr).
\end{eqnarray*}
Here, given $r\in\sR$ and $x\in\Rp$, $V_x^r(\cdot)$ is referred to
as a truncated
load process. Then, for $x\in\Rp$, $r\in\sR$ and $t\in[0,\infty)$,
\begin{eqnarray*}
\bigl\langle \chi1_{[0,x]},\sZhat^r(t)\bigr\rangle &=& \bigl
\langle \chi1_{[0,x]},\sZhat^r\bigl(\tau^r(t,x)
\bigr)\bigr\rangle + r \bigl(\Vbar_x^r(t)-
\Vbar_x^r\bigl(\tau^r(t,x)\bigr) \bigr)
\\
&&{}-r \theta^r(t,x),
\\
\bigl\langle \chi1_{[0,x]},\sZhat^r(t)\bigr\rangle &=& \bigl
\langle \chi1_{[0,x]},\sZhat^r\bigl(\tau^r(t,x)
\bigr)\bigr\rangle + \Vhat_x^r(t)-\Vhat_x^r
\bigl(\tau^r(t,x)\bigr)
\\
&&{}+ \bigl(\rho_x^r-1 \bigr) r \theta^r(t,x).
\end{eqnarray*}
Given $x\in\Rp$, $r\in\sR$ and $t\in[0,\infty)$, either $\tau
^r(t,x)=0$ or $\tau^r(t,x)>0$.
If the latter, then at time $\tau^r(t,x)$, either a job with total
processing time in $[0,x]$
arrives exogenously or a job with total processing time greater than
$x$ was in service
immediately before time $\tau^r(t,x)$, and its remaining processing
time at time $\tau^r(t,x)$
is $x$. Hence, for $x\in\Rp$, $r\in\sR$ and $t\in[0,\infty)$,
%
\begin{eqnarray}\label{eq:timetau}
\bigl\langle \chi1_{[0,x]},\sZhat^r\bigl(
\tau^r(t,x)\bigr)\bigr\rangle &\le& \bigl\langle \chi1_{[0,x]},
\sZhat^r(0)\bigr\rangle
\nonumber
\\
&&{}+ \frac{1}{r} \bigl(V_x^r\bigl(r^2
\tau^r(t,x)\bigr)-V_x^r\bigl(r^2
\tau^r(t,x)-\bigr) \bigr) + \frac{x}{r}
\\
&=& \bigl\langle \chi1_{[0,x]},\sZhat^r(0)\bigr\rangle
\nonumber
\\
&&{}+ r \bigl(\Vbar_x^r\bigl(\tau^r(t,x)
\bigr)-\Vbar_x^r\bigl(\tau^r(t,x)-\bigr)
\bigr) + \frac{x}{r}
\nonumber
\\
&=& \bigl\langle \chi1_{[0,x]},\sZhat^r(0)\bigr\rangle
\nonumber
\\
&&{}+ \Vhat_x^r\bigl(\tau^r(t,x)\bigr)-
\Vhat_x^r\bigl(\tau^r(t,x)-\bigr) +
\frac{x}{r}.
\nonumber
\end{eqnarray}
Therefore, for $x\in\Rp$, $r\in\sR$ and $t\in[0,\infty)$,
%
\begin{eqnarray}\label
{eq:RepWorkBelowX1}
\bigl\langle \chi1_{[0,x]},\sZhat^r(t)\bigr\rangle &\le&
\bigl\langle \chi1_{[0,x]},\sZhat^r(0)\bigr\rangle + r \bigl(
\Vbar_x^r(t)-\Vbar_x^r\bigl(
\tau^r(t,x)-\bigr) \bigr)
\nonumber
\\[-8pt]
\\[-8pt]
\nonumber
&&{} - r \theta^r(t,x) + \frac{x}{r},
\nonumber
\\
\label{eq:RepWorkBelowX2}
\bigl\langle \chi1_{[0,x]},\sZhat^r(t)\bigr\rangle &\le&
\bigl\langle \chi1_{[0,x]},\sZhat^r(0)\bigr\rangle +
\Vhat_x^r(t)-\Vhat_x^r\bigl(
\tau^r(t,x)-\bigr)
\nonumber
\\[-8pt]
\\[-8pt]
\nonumber
&&{}+ \bigl(\rho_x^r-1 \bigr) r \theta^r(t,x)
+ \frac{x}{r}.
\nonumber
\end{eqnarray}
We use \eqref{eq:RepWorkBelowX1} to prove the next lemma,
which specifies the asymptotic behavior of $\theta^r(\cdot,x)$ as
$r\to\infty$.
We use \eqref{eq:RepWorkBelowX2} to prove the subsequent lemma,
which specifies the asymptotic behavior of $\theta^r(\cdot
,l_{\varepsilon}^r)$ as
$r\to\infty$.

\begin{lemma}\label{lem:ThetaX}
For each $x\in\Rp$, as $r\to\infty$,
%
\begin{equation}
\label{eq:rlnr} c^r r\theta^r(\cdot,x)\Rightarrow0(\cdot).
\end{equation}
\end{lemma}

\begin{pf}
Given $x\in\Rp$, let $\rho_x=\lambda{\mathbb E}  [ v 1_{\{v
\le x\}}  ]$ and $\rho_x(t)=\rho_xt$ for all $[0,\infty)$.
Then, \eqref{eq:xTruncDiff} implies that, for each $x\in\Rp$, as
$r\to\infty$,
%
\begin{equation}
\label{eq:FLLNx} \Vbar_x^r(\cdot)\Rightarrow
\rho_x(\cdot).
\end{equation}
Fix $x\in\Rp$, $T>0$ and $\gamma>0$.
Note that $\rho_x<1$.
Let $\delta>0$ be such that $(1+\delta)\rho_x<1$. For $r\in\sR$, let
\begin{eqnarray*}
\Omega_0^r &=& \biggl\{ \bigl\langle
\chi1_{[0,x]},\sZtild^r (0)\bigr\rangle \le\frac
{\gamma}{2}
\biggr\},
\\
\Omega_1^r &=& \Bigl\{ \sup_{0\le s\le t\le T}
\Vbar_x^r(t)-\Vbar_x^r( s-)<(1+
\delta)\rho_x(t-s) \Bigr\},
\\
\Omega^r&=&\Omega_0^r\cap
\Omega_1^r.
\end{eqnarray*}
By \eqref{eq:initial2} and \eqref{eq:FLLNx},
\[
\lim_{r\to\infty}{\mathbb P} \bigl( \Omega^r \bigr)=1.
\]
By \eqref{eq:RepWorkBelowX1}, for each $r\in\sR$, on $\Omega^r$,
for each $t\in[0,T]$,
\[
\bigl\langle \chi1_{[0,x]},\sZhat^r(t)\bigr\rangle \le \bigl
\langle \chi1_{[0,x]},\sZhat^r (0)\bigr\rangle + \bigl((1+
\delta)\rho_x-1\bigr)r\theta^r(t,x)+\frac{x}{r}.
\]
But for each $r\in\sR$, $\langle \chi1_{[0,x]},\sZhat^r(t)\rangle
\ge0$
for each $t\in[0,T]$. Hence,
for each $r\in\sR$, on $\Omega^r$, for each $t\in[0,T]$,
\[
\bigl(1-(1+\delta)\rho_x\bigr)c^r r
\theta^r(t,x) \le \bigl\langle \chi1_{[0,x]},
\sZtild^r (0)\bigr\rangle + \frac{ c^r x}{r}.
\]
Recall that $S^{-1}(\cdot)$ is slowly varying so that $\lim_{y\to
\infty} S^{-1}(y) / y=0$. Hence $\lim_{r\to\infty} c^r/r=\lim_{r\to\infty} S^{-1}(r)/r=0$.
Then for $r\in\sR$ sufficiently large, on $\Omega^r$, for each $t\in[0,T]$,
\[
c^r r\theta^r(t,x)\le\frac{\gamma}{(1-(1+\delta)\rho_x)}.
\]
Since $\lim_{r\to\infty}{\mathbb P} (\Omega^r )=1$,
\eqref{eq:rlnr} holds.
\end{pf}

One feature of the SRPT discipline that is utilized in the above proof
is that by restricting
to jobs with remaining processing time in $[0,x]$ for a fixed $x$, the
workload process truncated
to jobs with remaining processing time in $[0,x]$ effectively behaves
as a subcritical queue.
We wish to obtain a version of Lemma~\ref{lem:ThetaX} on
$[0,l_{\varepsilon}^r]$
for fixed $\varepsilon>0$ with $r\to\infty$. Note that for
$\varepsilon>0$, $\lim_{r\to\infty} l_{\varepsilon}^r=\infty$.
Therefore, on such time intervals, the truncated workload process
approaches that of a critical queue.
This makes the verification of Lemma~\ref{lem:theta} a bit more
delicate, and the rate of convergence result obtained
is not as rapid. For this, for $\varepsilon>0$, $r\in\sR$ and $t\in
[0,\infty)$, we adopt the shorthand notation
\[
\tau_{\varepsilon}^r(t)=\tau^r\bigl(t,l_{\varepsilon}^r
\bigr) \quad\mbox{and}\quad \theta_{\varepsilon}^r(t)=\theta^r
\bigl(t,l_{\varepsilon}^r\bigr).
\]

\begin{lemma}\label{lem:theta}
For $\varepsilon>0$, as $r\to\infty$,
\[
\bigl(c^r \bigr)^{2+\varepsilon} \theta_{\varepsilon}^r(
\cdot )\Rightarrow0(\cdot).
\]
\end{lemma}

\begin{pf} Fix $\varepsilon>0$ and $t\in[0,\infty)$.
Given $r\in\sR$, we take $x=l_{\varepsilon}^r$ in \eqref
{eq:RepWorkBelowX2}, and then we
subtract and add $\rho^r r\theta_{\varepsilon}^r(t)$, and use \eqref
{eq:rhominusrhox} and
the fact that $S(l_{\varepsilon}^r)=r(c^r)^{-2-\varepsilon} $ to
obtain that for $r\in\sR$,
%
\begin{eqnarray}\label{eq:RepWork}
\bigl\langle \chi1_{[0,l_{\varepsilon}^r]},\sZhat^r(t)\bigr\rangle &\le&
\bigl\langle \chi1_{[0,l_{\varepsilon}^r]},\sZhat^r(0)\bigr\rangle +
\Vhat_{l_{\varepsilon}^r}^r(t)-\Vhat_{l_{\varepsilon}^r}^r\bigl(\tau
_{\varepsilon}^r(t)-\bigr)
\nonumber
\\[-8pt]
\\[-8pt]
\nonumber
&&{}- \lambda^r \bigl(c^r\bigr)^{2+\varepsilon}
\theta_{\varepsilon}^r(t) + \bigl(\rho^r-1 \bigr) r
\theta_{\varepsilon}^r(t) + \frac{ l_{\varepsilon}^r }{r}.
\nonumber
\end{eqnarray}
We have that $\langle \chi1_{[0,l_{\varepsilon}^r]},\sZhat
^r(t)\rangle \ge
0$ and $\theta_{\varepsilon}^r(t)\ge0$
for all $r\in\sR$. This together with the fact that $l_{\varepsilon
}^r<c^r$ implies that, for all $r\in\sR$,
%
\begin{eqnarray}\quad
\label{eq:KeyINEQ1} 0&\le&\lambda^r \bigl(c^r
\bigr)^{2+\varepsilon} \theta_{\varepsilon}^r(t)
\nonumber
\\[-8pt]
\\[-8pt]
\nonumber
&\le& \bigl\langle
\chi1_{[0,l_{\varepsilon}^r]},\sZhat^r(0)\bigr\rangle + \Vhat_{l_{\varepsilon}^r}^r(t)-
\Vhat_{l_{\varepsilon}^r}^r\bigl(\tau _{\varepsilon}^r(t)-
\bigr)+ \bigl(\rho^r-1 \bigr) r \theta_{\varepsilon}^r(t) +
\frac{c^r}{r}.
\end{eqnarray}
%
Upon dividing by $(c^r)^{2+\varepsilon}$ and using $\lim_{r\to\infty
}c^r=\infty$, \eqref{eq:param}, \eqref{eq:initial2} and \eqref
{eq:TruncDiff},
we see that, as $r\to\infty$,
%
\begin{equation}
\label{eq:step1} \theta_{\varepsilon}^r(\cdot)\Rightarrow0(\cdot).
\end{equation}
Hence, by \eqref{eq:TruncDiff} and the fact that $V^*(\cdot)$ is
continuous, as $r\to\infty$,
%
\begin{equation}
\label{eq:DiffZero} \Vhat_{l_{\varepsilon}^r}^r(\cdot)-\Vhat_{l_{\varepsilon
}^r}^r
\bigl(\tau_{\varepsilon}^r(\cdot)-\bigr)\Rightarrow0(\cdot).
\end{equation}
Then letting $r\to\infty$ in \eqref{eq:KeyINEQ1} and
using \eqref{eq:param}, \eqref{eq:initial2}, \eqref{eq:step1},
\eqref{eq:DiffZero} and the fact that $c^r=S^{-1}(r)$ and
$S^{-1}(\cdot)$
is slowly varying completes the proof.
\end{pf}

\subsubsection{Truncated queue length process asymptotics}\label{sec:SRPTQ}

We are prepared to use Lemma~\ref{lem:ThetaX} to verify that the total
mass in a fixed width neighborhood of the origin vanishes as $r$ tends
to infinity.

\begin{lemma}\label{lem:Queuex}
For all $x\in\Rp$, as $r\to\infty$,
%
\begin{equation}
\label{eq:NearOriginQ} \bigl\langle 1_{[0,x]},\sZtild^r(\cdot)\bigr
\rangle \Rightarrow0(\cdot).
\end{equation}
\end{lemma}

\begin{pf} Fix $x\in\Rp$ and $T>0$. By ignoring any processing that occurs
in $(r^2\tau^r(t,x),r^2t]$, for $r\in\sR$ and $t\in[0,\infty)$, we
have that
\[
\bigl\langle 1_{[0,x]},\sZtild^r(t)\bigr\rangle \le \bigl
\langle 1_{[0,x]},\sZtild^r\bigl(\tau^r(t,x)\bigr)
\bigr\rangle + c^r r \bigl( \Ebar^r(t)-\Ebar^r
\bigl(\tau^r(t,x)\bigr) \bigr).
\]
Further, by using arguments similar to those that yielded \eqref
{eq:timetau}, for $r\in\sR$ and $t\in[0,\infty)$,
\[
\bigl\langle 1_{[0,x]},\sZtild^r\bigl(\tau^r(t,x)
\bigr)\bigr\rangle \le \bigl\langle 1_{[0,x]},\sZtild^r(0)\bigr
\rangle + c^rr \bigl(\Ebar^r\bigl(\tau^r(t,x)
\bigr)-\Ebar^r\bigl(\tau^r(t,x)-\bigr) \bigr) +
\frac{c^r}{r}.
\]
Then, for $r\in\sR$ and $t\in[0,\infty)$, we have that
%
\begin{equation}\quad
\label{eq:QueueRep} \bigl\langle 1_{[0,x]},\sZtild^r(t)\bigr\rangle
\le \bigl\langle 1_{[0,x]},\sZtild^r(0)\bigr\rangle +
c^r r \bigl( \Ebar^r(t)-\Ebar^r\bigl(
\tau^r(t,x)-\bigr) \bigr) +\frac{c^r}{r}.
\end{equation}
Fix $\gamma>0$. For $r\in\sR$, let
\begin{eqnarray*}
\Omega_0^r &=& \biggl\{ \bigl\langle 1_{[0,x]},
\sZtild^r (0)\bigr\rangle \le\frac{\gamma
}{3} \biggr\},
\\
\Omega_1^r &=& \Bigl\{ \sup_{0\le s\le t\le T}
\bigl(\Ebar^r(t)-\Ebar^r( s-) \bigr)<2\lambda(t-s) \Bigr
\},
\\
\Omega_2^r&=& \biggl\{ \sup_{t\in[0,T]}
\theta^r(t,x)<\frac{\gamma
}{6\lambda c^r r} \biggr\},
\\
\Omega^r&=&\Omega_0^r\cap
\Omega_1^r\cap\Omega_2^r.
\end{eqnarray*}
By \eqref{eq:FLLNE}, \eqref{eq:initial2} and \eqref{eq:rlnr}, $\lim_{r\to\infty}{\mathbb P} ( \Omega^r )=1$.
Then since $c^r=S^{-1}(r)$ and $S^{-1}(\cdot)$ is slowly varying, it
follows that, on $\Omega^r$, for $r$ sufficiently large,
\[
\sup_{t\in[0,T]}\bigl\langle 1_{[0,x]},\sZtild^r(t)
\bigr\rangle \le \gamma.
\]
Since $\gamma>0$ was arbitrary, the proof is complete.
\end{pf}


\subsubsection{Truncated workload process asymptotics}\label{sec:SRPTWBelow}

We are prepared to use Lemma~\ref{lem:theta} to verify that the total
work in a growing neighborhood of the origin vanishes
as $r$ tends to infinity.

\begin{lemma}\label{lem:Worktild}
For all $\varepsilon>0$, as $r\to\infty$,
\[
\bigl\langle \chi1_{[0,l_{\varepsilon}^r]},\sZtild^r(\cdot)\bigr\rangle
\Rightarrow 0(\cdot).
\]
\end{lemma}

Before proving Lemma~\ref{lem:Worktild}, we begin with an observation.
By \eqref{eq:RepWork}, for each $\varepsilon>0$, $r\in\sR$ and
$t\in[0,\infty)$,
%
\begin{eqnarray}\label{eq:RepWork2}
\bigl\langle \chi1_{[0,l_{\varepsilon}^r]},\sZtild^r(t)\bigr\rangle &\le&
\bigl\langle \chi1_{[0,l_{\varepsilon}^r]},\sZtild^r(0)\bigr\rangle +
c^r \bigl(\Vhat_{l_{\varepsilon}^r}^r(t)-\Vhat_{l_{\varepsilon
}^r}^r
\bigl(\tau_{\varepsilon}^r(t)-\bigr) \bigr)
\\
&&{}- \lambda^r \bigl(c^r\bigr)^{3+\varepsilon}
\theta_{\varepsilon}^r(t) + c^r \bigl(
\rho^r-1 \bigr) r \theta_{\varepsilon}^r(t) +
\frac{c^r l_{\varepsilon}^r} {r}.
\nonumber
\end{eqnarray}
We argue that each term on the right-hand side converges in
distribution to the zero process.
We begin by proving the following lemma.

\begin{lemma}\label{lem:Holder}
For each $\varepsilon>0$, as $r\to\infty$,
\[
c^r \bigl( \Vhat_{l_{\varepsilon}^r}^r(\cdot)-\Vhat
_{l_{\varepsilon}^r}^r\bigl(\tau_{\varepsilon}^r(\cdot)-
\bigr) \bigr)\Rightarrow0(\cdot).
\]
\end{lemma}

\begin{pf} Fix $T,\varepsilon>0$.
Recall that Brownian motion is H\"older continuous with
exponent $\gamma$ for any $0<\gamma<1/2$.
Fix $0<\gamma<1/2$ such that $\gamma(2+\varepsilon)>1$.
For $M\in\N$, let
\[
\Omega(M)= \bigl\{ \bigl\llvert V^*(t)-V^*(s-) \bigr\rrvert <
M(t-s)^{\gamma
}\mbox{ for all }0\le s\le t\le T \bigr\}.
\]
We have that $\Omega(M)\subset\Omega(M+1)$ for all $M\in\N$ and
\[
{\mathbb P} \biggl( \bigcup_{M\in\N} \Omega(M)
\biggr)=1.
\]
Hence given $\eta>0$, there exists $M_{\eta}\in\N$ such that
\[
{\mathbb P} \bigl( \Omega(M_{\eta}) \bigr)\ge1-\eta.
\]
Given $r\in\sR$ and $M\in\N$, let
\[
\Omega^r(M)= \bigl\{ \bigl\llvert \Vhat_{l_{\varepsilon}^r}^r(t)-
\Vhat _{l_{\varepsilon}^r}^r(s-) \bigr\rrvert < M(t-s)^{\gamma}\mbox{
for all }0\le s\le t\le T \bigr\}.
\]
For each $M\in\N$, the set $A(M)$, given by
\[
A(M)= \bigl\{ f\in\bD\bigl([0,T],\R\bigr) \dvtx \bigl\llvert f(t)-f(s-) \bigr
\rrvert < M(t-s)^{\gamma}\mbox{ for all }0\le s\le t\le T \bigr\},
\]
is open in the uniform topology. Hence, 
\eqref{eq:TruncDiff} and the Portmanteau theorem imply that
\[
\liminf_{r\to\infty}{\mathbb P} \bigl(\Omega^r(M_{\eta})
\bigr)\ge{\mathbb P} \bigl(\Omega(M_{\eta}) \bigr) \ge1-\eta.
\]
For $r\in\sR$, let
\[
\Omega_1^r= \biggl\{ \sup_{t\in[0,T]}
\theta_{\varepsilon
}^r(t)<\frac{1}{\sqrt[\gamma]{M_{\eta}} (c^r)^{2+\varepsilon}} \biggr\}.
\]
By Lemma~\ref{lem:theta},
\[
\lim_{r\to\infty}{\mathbb P} \bigl( \Omega_1^r
\bigr)=1.
\]
Then
\[
\liminf_{r\to\infty}{\mathbb P} \bigl(\Omega^r(M_{\eta})
\cap \Omega_1^r \bigr)\ge1-\eta.
\]
Given $r\in\sR$, set
\[
\Omega_2^r = \biggl\{ \sup_{t\in[0,T]}
c^r \bigl\llvert \Vhat^r(t)-\Vhat^r\bigl(\tau
_{\varepsilon}^r(t)-\bigr)\bigr\rrvert < \frac{ c^r }{ (c^r)^{(2+\varepsilon)\gamma}}=
\frac{ 1 }{
(c^r)^{(2+\varepsilon)\gamma-1}} \biggr\}.
\]
Then, for $r\in\sR$,
\[
\Omega^r(M_{\eta})\cap\Omega_1^r
\subset\Omega_2^r.
\]
Hence
\[
\liminf_{r\to\infty}{\mathbb P} \bigl(\Omega_2^r
\bigr)\ge1-\eta.
\]
But, for $r\in\sR$, $\Omega_2^r$ does not depend on $\eta$.
Therefore, we may let $\eta$ decrease to zero
so that
\[
\liminf_{r\to\infty}{\mathbb P} \bigl(\Omega_2^r
\bigr)= 1.
\]
Fix $\delta>0$. Given $r\in\sR$, let
\[
\Omega_3^r= \Bigl\{ \sup_{t\in[0,T]}
c^r \bigl\llvert \Vhat^r(t)-\Vhat ^r\bigl(
\tau_{\varepsilon}^r(t)-\bigr)\bigr\rrvert <\delta \Bigr\}.
\]
Since $\lim_{r\to\infty}c^r=\infty$ and $(2+\varepsilon)\gamma
-1>0$, it follows that for $r$ sufficiently large
$\Omega_2^r\subset\Omega_3^r$. Therefore, $\liminf_{r\to\infty
}{\mathbb P} (\Omega_3^r )= 1$.
Since $T,\varepsilon,\delta>0$ were arbitrary, Lemma~\ref
{lem:Holder} holds.
\end{pf}

\begin{corollary}\label{cor:3rate} For each $\varepsilon>0$, as $r\to
\infty$,
\[
\lambda^r \bigl(c^r\bigr)^{3+\varepsilon}
\theta_{\varepsilon}^r(\cdot )\Rightarrow0(\cdot).
\]
\end{corollary}

\begin{pf} Fix $\varepsilon>0$. By \eqref{eq:RepWork2}, the fact that
$\langle \chi1_{[0,l_{\varepsilon}^r]},\sZtild^r(t)\rangle \ge0$ and
$\theta_{\varepsilon}^r(t)\ge0$
for all $r\in\sR$ and $t\in[0,\infty)$ and $l_{\varepsilon}^r< c^r$
for all $r\in\sR$,
we have that, for all $r\in\sR$ and $t\in[0,\infty)$,
%
\begin{eqnarray}
0&\le& \lambda^r \bigl(c^r\bigr)^{3+\varepsilon}
\theta_{\varepsilon}^r(t)\nonumber
\\
&\le& \bigl\langle \chi1_{[0,l_{\varepsilon}^r]},\sZtild^r(0)\bigr\rangle +
c^r \bigl(\Vhat_{l_{\varepsilon}^r}^r(t)-\Vhat_{l_{\varepsilon
}^r}^r
\bigl(\tau_{\varepsilon}^r(t)-\bigr) \bigr)
\\
&&{}+ c^r \bigl(\rho^r-1 \bigr) r \theta_{\varepsilon}^r(t)
+ \frac{  (c^r )^2 }{r}.
\nonumber
\end{eqnarray}
The result follows from this, \eqref{eq:initial2}, Lemma~\ref
{lem:Holder}, \eqref{eq:param}, Lemma~\ref{lem:theta} and the fact
that $c^r=S^{-1}(r)$ and $S^{-1}(\cdot)$ is slowly varying.
\end{pf}

\begin{pf*}{Proof of Lemma~\ref{lem:Worktild}}
Fix $\varepsilon>0$.
The result follows
by combining \eqref{eq:RepWork2}, \eqref{eq:initial2}, Lemma~\ref
{lem:Holder}, Corollary~\ref{cor:3rate},
\eqref{eq:param}, Lemma~\ref{lem:theta}, $l_{\varepsilon}^r<c^r$ for
$r\in\sR$,
$c^r=S^{-1}(r)$ for $r\in\sR$ and $S^{-1}(\cdot)$ is slowly varying.
\end{pf*}

\begin{appendix}\label{app}

\section*{Appendix: Behavior of truncated load processes}
The following result is well known and follows from \cite{ref:P}, Theorem~3.1, used to extend \cite{ref:B}, Section~17.3.

\begin{propositionn}\label{prop:basic}
For each $r\in\sR$, let $\{x_k^r\}_{k=1}^{\infty}$ be an independent and
identically distributed sequence of nonnegative random variables
with finite mean $m^r$ and finite standard deviation $\sigma^r$ that is
independent of $E^r(\cdot)$. Suppose that for some finite nonnegative
constants $m$ and $\sigma$, $\lim_{r\to\infty}m^r=m$ and $\lim_{r\to\infty}\sigma^r=\sigma$.
Further assume that for each $\delta>0$,
\[
\lim_{r\to\infty}\E \bigl[ \bigl(x_1^r-m^r
\bigr)^2 \bigl\llvert x_1^r-m^r
\bigr\rrvert > r\delta \bigr]=0.
\]
For $r\in\sR$, $n\in\N$ and $t\in[0,\infty)$, let
\[
X^r(n)=\sum_{k=1}^n
x_k^r\quad \mbox{and}\quad \Xhat^r(t)=
\frac{X^r(\lfloor r^2t\rfloor) -\lfloor r^2t\rfloor m^r}{r}.
\]
Then, as $r\to\infty$, $(\Ehat^r(\cdot),\Xhat^r(\cdot
))\Rightarrow(E^*(\cdot),X^*(\cdot))$,
where $E^*(\cdot)$ is given by \eqref{eq:FCLTE}, and $X^*(\cdot)$ is
a Brownian motion starting from zero with zero drift and
variance $\sigma^2$ per unit time, that is independent of $E^*(\cdot)$.
Furthermore, as $r\to\infty$,
\[
\frac{X^r(r^2\Ebar^r(\cdot))-r^2\lambda^r(\cdot)m^r}{r} \Rightarrow X^*\bigl(\lambda(\cdot)\bigr)+m E^*(\cdot),
\]
where for each $r\in\sR$ and $t\in[0,\infty)$, $\lambda
^r(t)=\lambda
^r t$ and $\lambda(t)=\lambda t$.
\end{propositionn}

Recall that, for $r\in\sR$ and $x\in\Rp$,
\[
\Vhat_x ^r(\cdot) = \frac{\sum_{i=1}^{r^2\Ebar^r(\cdot)} v_i1_{\{v_i\le x\}}
-
r^2\lambda^r(\cdot) {\mathbb E} [ v 1_{\{v\le x\}} ] }{r}.
\]
Proposition~\ref{prop:basic} implies that for each $x\in\Rp$,
as $r\to\infty$,
%
\setcounter{equation}{0}
\begin{equation}
\label{eq:xTruncDiff} \Vhat_x^r(\cdot)\Rightarrow
V_x^*(\cdot),
\end{equation}
where $V_x^*(\cdot)$ is a Brownian motion starting from zero with
drift zero
and finite variance per unit time.
Similarly, Proposition~\ref{prop:basic} together with
$0 \le\E[v 1_{\{v\le l_{\varepsilon}^r\}}]\le\E[v] $ and $0 \le\E
[v^2 1_{\{v\le l_{\varepsilon}^r\}}]\le\E[v^2] $ for all
$\varepsilon>0$ and $r\in\sR$ and the monotone convergence theorem
implies that for each $\varepsilon>0$, as $r\to\infty$,
%
\begin{equation}
\label{eq:TruncDiff} \Vhat_{l_{\varepsilon}^r}^r(\cdot)\Rightarrow V^*(\cdot),
\end{equation}
where $V^*(\cdot)$ is a Brownian motion starting from zero with drift zero
and finite variance per unit time.
\end{appendix}

\section*{Acknowledgment} The author would like to thank ViaSat Inc.
for generously funding
undergraduate research assistants Richard Hunperger and Sean Malter who
developed
code and performed simulations that helped the author formulate the
statement of
Theorem~\ref{thrm:main}.


%





\printaddresses
\end{document}